\documentclass[11pt]{article}

\usepackage{amssymb,amsmath,euscript,bbm,xcolor,graphicx,epstopdf}
\RequirePackage[numbers]{natbib}

\usepackage{amsmath,amssymb,amsthm,bm,euscript,bbm}
\usepackage{verbatim}
\usepackage{graphicx}
\usepackage{epstopdf}
\usepackage{color}
\usepackage{pstricks,pst-node,pst-plot}
\usepackage{rotating}
\usepackage{enumitem}
\newtheorem{proposition}{Proposition}[section]
\newtheorem{lemma}[proposition]{Lemma}
\newtheorem{theorem}[proposition]{Theorem}

\newtheorem{corollary}[proposition]{Corollary}

\def\la{\lambda}

\def\R{{\mathbb R}}

\def\E{{\mathbb E}}
\def\P{{\mathbb P}}

\makeatletter \@addtoreset{equation}{section} \makeatother
\newenvironment{remark}{%
	\vspace{0.3cm} \pagebreak [2]%
	\par%
	\refstepcounter{proposition}
	\noindent%
	{\bf Remark~\theproposition\  }}{\qed}%
\begin{document}
	
	\title {On local maxima of smooth Gaussian nonstationary processes and stationary planar fields with trends}
	\author{Dan Cheng\\ Arizona State University}
	
	\date{}
	
	\maketitle
	
	\begin{abstract}
		We present exact formulas for both the expected number and the height distribution of local maxima (peaks) in two distinct categories of smooth, non-centered Gaussian fields: (i) nonstationary Gaussian processes and (ii) stationary planar Gaussian fields. For case (i), we introduce a novel parameter related to conditional correlation that significantly simplifies the computation of these formulas. Notably, the peak height distribution is solely dependent on this single parameter. In case (ii), traditional methods involving GOE random matrices are ineffective for non-isotropic fields with mean functions. To address this, we apply specific transformations that enable the derivation of formulas using generalized chi-squared density functions. These derived results provide essential tools for calculating p-values and power in applications of signal  and change point detection within environments characterized by non-isotropic Gaussian noise.
	\end{abstract}
	
	
	
	\section{Introduction}
	
	The study of local maxima (peaks) of Gaussian random fields, particularly their expected number and height distribution, is a significant topic in probability theory \cite{AT07, CL67}, with extensive applications across various fields such as statistics \cite{CS17, ChangePoint, CS19sphere}, physics \cite{Annibale:2003}, neuroimaging \cite{Taylor:2007, Pantazis:2005}, oceanography \cite{Lindgren82}, and astronomy \cite{AstrophysJ85}. Researchers from these diverse areas have developed powerful tools to address this problem, notably the well-known Kac-Rice formula \cite{AT07}. While the Kac-Rice formula provides an implicit solution for the expected number of local maxima, explicit evaluation remains challenging due to complex computations involving Hessian matrices. An exception exists for centered isotropic Gaussian fields, where techniques from random matrix theory (specifically GOE and GOI matrices) facilitate the handling of Hessian computations, thus yielding explicit formulae \cite{CS15, CS18}. This paper extends the investigation to nonstationary Gaussian processes in 1D and stationary Gaussian fields in 2D with mean functions, deriving exact formulae for the expected number and height distribution of local maxima for a broader class of Gaussian fields. These derived results will provide essential tools for calculating p-values and power in signal and change point detection under non-isotropic Gaussian noise.
	
	Specifically, let $X=\{X(t), t\in T\}$ be a noncentered, smooth and unit-variance Gaussian random field living on an \textit{open and bounded} parameter set $T$ in $\R$ or $\R^2$. Let
	\begin{equation}\label{Eq:M-Mu}
		\begin{split}
		M(X, T) &= \# \left\{ t\in T: \nabla X(t)=0, \nabla^2 X(t)\prec 0 \right\}, \\
		M_u(X, T) &= \# \left\{ t\in T: X(t)\geq u, \nabla X(t)=0, \nabla^2 X(t)\prec 0 \right\}, 
		\end{split}
	\end{equation}
	where $\nabla X(t)$ and $\nabla^2 X(t)$ are the gradient and Hessian of $X$, respectively, and $\nabla^2 X(t)\prec 0$ indicates that the Hessian is negative definite. Thus, $M(X, T)$ and $M_u(X, T)$ represent the number of local maxima and the number of local maxima exceeding $u$ of $X$ over $T$, respectively. Calculating $\E[M(X, T)]$ and $\E[M_u(X, T)]$ is often complicated due to the expectations involving $\nabla^2 X(t)\prec 0$, and this complexity increases significantly when incorporating a mean function. Thus, existing results typically focus on centered and isotropic Gaussian fields. In this paper, we address 1D nonstationary and 2D stationary Gaussian fields, presenting methods to simplify the calculations involving the Hessian. Our main contributions are outlined as follows.
	
	For the 1D Gaussian process, the absence of stationarity might suggest that the formula for the expected number of local maxima would be complex. However, by introducing a new parameter $\rho_t$ (cf. \eqref{eq:rho_1} and \eqref{eq:rho} below) based on conditional correlation, we significantly simplify the derived formulas; see Theorem \ref{thm:1D}. Importantly, this new parameter generalizes the $\kappa$ parameter used for isotropic Gaussian fields \cite{CS15,CS18} and provides a meaningful interpretation for the $\kappa$ parameter as well.
	
	For the 2D Gaussian fields, the general formula for the stationary nonisotropic case remains unknown. We leverage the fact that $\nabla^2 X(t)\prec 0$ is equivalent to both diagonal entries being negative and the determinant being positive. By transforming the distribution of the determinant of $\nabla^2 X(t)$ into a generalized chi-squared distribution, we are able to handle the expectation involving $\nabla^2 X(t)$. The derived formulas in Theorem \ref{theorem:2D} are expressed in terms of the generalized chi-squared density. This technique also allows us to derive the formula in Corollary \ref{cor:Isotropic} for isotropic planar Gaussian fields with a mean function. Furthermore, we verify that the derived formula using the chi-squared distribution aligns with the formula derived using the random matrix technique in \cite{CS18}.
	
	  The peak height distribution of $X$ at a point $t$ is defined as the probability that the peak height exceeds a fixed threshold at $t$, conditioned on the event that $t$ is a local maximum of $X$. More rigorously, it is defined as
	\begin{equation*}
		F_t(u) = \P[X(t)>u \ | \ t \text{ is a local maximum of } X ].
	\end{equation*}
	It has been shown in \cite{CS15} that the peak height distribution of $X$ at $t$ is given by
	\begin{equation}\label{eq:Palm distr Euclidean}
		F_t(u) = \frac{\E\left[|{\rm det} (\nabla^2 X(t))|\mathbbm{1}_{\{X(t)> u\}} \mathbbm{1}_{\{\nabla^2 X(t)\prec 0\}}|\nabla X(t)=0\right]}{\E\left[|{\rm det} (\nabla^2 X(t))|\mathbbm{1}_{\{\nabla^2 X(t)\prec 0\}} | \nabla X(t)=0\right]}, \quad  u\in \R.
	\end{equation}
	Taking the derivative, we obtain the peak height density, denoted by $h_t(x) = -F_t'(x)$, $x\in \R$. 
	
	Similar to the computation of $\E[M(X, T)]$ and $\E[M_u(X, T)]$, it can be seen from \eqref{eq:Palm distr Euclidean} that the expectation involving the determinant and indicator function on the Hessian is challenging. Here, for the 1D case, we have established in Corollaries \ref{cor:1D_PHD} and \ref{cor:1D-1} that the peak height distribution for a centered nonstationary Gaussian process depends solely on the parameter $\rho_t$. This result is particularly valuable in statistical applications, as it enables inference on the peak height distribution based on the estimation of $\rho_t$ only. For 2D stationary Gaussian fields, we employ generalized chi-squared densities to derive the peak height density in Theorem \ref{theorem:2D}.
	
	Throughout this paper, we assume that the mean function $m(t)=\E[X(t)]$ is twice differentiable, and that the field $X$ satisfies the regularity and  smoothness conditions specified in Corollary 11.3.2 in \cite{AT07}. Roughly speaking, it requires that the joint distributions of the gradient and Hessian of $X$ are nondegenerate and that the sample paths of $X$ are in $C^2$ almost surely; see more details in Sections \ref{sec:1D} and \ref{sec:2D} below. These conditions imply that $X$ is almost surely a Morse function, allowing us to apply the Kac-Rice formula to compute the expected number of local maxima. Denote by $\phi(x)$ and $\Phi(x)$ the pdf and cdf of the standard normal distribution $N(0,1)$, respectively. Let $\Psi(x)=1-\Phi(x)$ be tail probability of $N(0,1)$. Furthermore, let $Y_1\overset{d}{=}Y_2$ denote that random variables or vectors $Y_1$ and $Y_2$ have the same distribution.
	
	\section{Nonstationary Gaussian Processes} \label{sec:1D}
	In this section, we study the case that $\{X(t), t\in T\subset \R\}$ is a smooth, unit-variance, nonstationary 1D Gaussian process. Specifically, following Corollary 11.3.2 in \cite{AT07}, we assume that the distributions of $(X'(t), X''(t))$ are nondegenerate for all $t\in T$, and that there exist constants $K, \alpha>0$ such that for all $t,s\in T$,
	\[
	|{\rm Var}(X''(t)) + {\rm Var}(X''(s))-2{\rm Cov}(X''(t), X''(s))| \le K|\log |t-s||^{-(1+\alpha)}.
	\]
	 Suppose $X$ has the mean function $\E[X(t)]=m(t)$. Let 
	\begin{equation}\label{eq:delta_1}
		\la_1(t)={\rm Var}(X'(t)), \quad \la_2(t)={\rm Var}(X''(t)), \quad   \delta_t^2 = {\rm Var}(X''(t) | X'(t)=0).
	\end{equation}
Let $\rho_t$ be the \textit{conditional correlation} of $X(t)$ and $X''(t)$ given $X'(t)=0$, that is,
\begin{equation}\label{eq:rho_1}
		\rho_t= {\rm Corr}(X(t), X''(t) | X'(t)=0).
\end{equation} 
Note that, if $X$ is a stationary Gaussian process, then $X'(t)$ is independent of $X(t)$ and $X''(t)$ for each $t$, and hence $\delta_t^2 = {\rm Var}(X''(t))=\la_2(t)$ and $\rho_t={\rm Corr}(X(t), X''(t))=-\la_1(t)/\sqrt{\la_2(t)}$, where the label $t$ can be omitted since the variance and correlation do not depend on $t$ under stationarity. For a general nonstationary Gaussian process, we have the following result.
\begin{lemma}
	Let $\{X(t), t\in T\subset \R\}$ be a smooth, unit-variance, nonstationary Gaussian process. Then 
	\begin{equation}\label{eq:delta}
		\delta_t^2 = \la_2(t) - \la_1'(t)^2/(4\lambda_1(t))
	\end{equation}
and 
\begin{equation}\label{eq:rho}
	\rho_t= \frac{-\la_1(t)}{\sqrt{\la_2(t) - \la_1'(t)^2/(4\lambda_1(t))}} \in [-1, 0).
\end{equation} 
\end{lemma}
\begin{proof}
	Note that ${\rm Var}(X(t)) = \E[(X(t)-m(t))^2] \equiv 1$. Taking the derivative on both sides, we obtain $\E[(X(t)-m(t))(X'(t)-m'(t))]=0$; that is, ${\rm Cov}(X(t), X'(t))=0$. Taking the derivative again yields
	\[
	\E[(X'(t)-m'(t))(X'(t)-m'(t))] + \E[(X(t)-m(t))(X''(t)-m''(t))]=0,
	\]
	implying ${\rm Cov}(X(t), X''(t))=-{\rm Var}(X'(t))=-\la_1(t)$. Similarly, taking the derivative on both sides of ${\rm Var}(X'(t))=\la_1(t)$, we obtain ${\rm Cov}(X'(t), X''(t))=\la_1'(t)/2$. Thus, the variance-covariance matrix of the Gaussian vector $(X(t), X'(t), X''(t))$ is given by
	\begin{equation}\label{eq:var-cov}
		\begin{split}
			\begin{pmatrix}
				1 & 0 & -\la_1(t) \\
				0 & \la_1(t) & \la_1'(t)/2 \\
				-\la_1(t) & \la_1'(t)/2 & \la_2(t)
			\end{pmatrix}
		\end{split}.
	\end{equation}  
	Applying the well-known formula for the conditional distribution of a  multivariate normal distribution, we have the following conditional variances and covariance:
	\begin{equation}\label{eq:delta_2}
		\begin{split}
			{\rm Var}(X(t) | X'(t)=0) &= 1,\\
			{\rm Var}(X''(t) | X'(t)=0) &= \la_2(t) - \la_1'(t)^2/(4\lambda_1(t)),\\
			{\rm Cov}(X(t), X''(t) | X'(t)=0) &= -\lambda_1(t),
		\end{split}
	\end{equation} 
	which implies \eqref{eq:delta} and \eqref{eq:rho}. Finally, as a variance function, $\la_1(t)>0$, thus $\rho_t<0$.
\end{proof}
	
	\subsection{The nondegenerate case: $-1<\rho_t<0$}
	We first consider the nondegenerate case where $\rho_t\neq -1$. The following result provides the exact formulas for  computing $\E[M(X, T)]$ and $\E[M_u(X, T)]$. 
	\begin{theorem}\label{thm:1D}
		Let $\{X(t), t\in T\subset \R\}$ be a smooth, unit-variance, nonstationary Gaussian process with mean function $m(t)$. Suppose $-1<\rho_t<0$ for all $t\in T$. Then the expected numbers of local maxima are given by
		\begin{align}		
				\E[M(X, T)]&=\int_T \left[\phi(b_t)+b_t\Phi(b_t)\right]\frac{\delta_t}{\sqrt{2\pi \la_1(t)}} e^{-\frac{m'(t)^2}{2\la_1(t)}}dt, \label{eq:LocalMax-1D-1}	\\
				\E[M_u(X, T)]&= \int_T \Bigg[ \phi(b_t)\Psi\left(\frac{u-m(t)-\rho_tb_t}{\sqrt{1-\rho_t^2}}\right) - \rho_t\phi(u-m(t))\Phi\left(\frac{b_t-\rho_t(u-m(t))}{\sqrt{1-\rho_t^2}}\right) \label{eq:LocalMax-1D-2} \\
				&\qquad \quad  + b_t \P(\xi_1>u-m(t), \xi_2 <b_t)\Bigg]\frac{\delta_t}{\sqrt{2\pi \la_1(t)}} e^{-\frac{m'(t)^2}{2\la_1(t)}}dt, \nonumber
		\end{align} 
	where $\delta_t$ and $\rho_t$ are given in \eqref{eq:delta} and \eqref{eq:rho}, and
		\begin{equation}\label{eq:xi}
			\begin{split}
				b_t = \frac{-m''(t)+m'(t)\la_1'(t)/(2\lambda_1(t))}{\delta_t}, \quad   (\xi_1, \xi_2) \sim N\left(\begin{pmatrix}
					0 \\
					0
				\end{pmatrix},\, \begin{pmatrix}
					1 & \rho_t\\
					\rho_t & 1
				\end{pmatrix}\right).
			\end{split}
		\end{equation} 
	\end{theorem}
	\begin{proof}  
		By \eqref{eq:var-cov} and \eqref{eq:delta_2}, we have the following joint conditional distribution:
		\begin{equation}\label{eq:cov-conditional-m}
			\begin{split}
				(X(t), X''(t) |  X'(t)=0) \sim N\left(\begin{pmatrix}
					m(t) \\
					-\delta_tb_t
				\end{pmatrix},\, \begin{pmatrix}
					1 & -\lambda_1(t)\\
					-\lambda_1(t) & \delta_t^2 
				\end{pmatrix}\right),
			\end{split}
		\end{equation} 
		which implies $(X''(t) |  X'(t)=0) \sim N(-\delta_tb_t,\, \delta_t^2)$. Therefore,
		\begin{equation}\label{eq:J1-1D}
			\begin{split}
				\E\left[|X''(t)|\mathbbm{1}_{\{X''(t) < 0\}} | X'(t)=0\right] 
				= \int_{-\infty}^0 \frac{-x}{\sqrt{2\pi}\delta_t}e^{-\frac{(x+\delta_tb_t)^2}{2\delta_t^2}}dx =  \delta_t\left[\phi(b_t)+b_t\Phi(b_t)\right].
			\end{split}
		\end{equation} 
		Note that $p_{X'(t)}(0) = [2\pi \la_1(t)]^{-1/2} e^{-m'(t)^2/\la_1(t)}$. By the Kac-Rice formula (cf.  Theorem 11.2.1 in \cite{AT07}),
		\begin{equation*}
			\begin{split}
				\E[M(X, T)]&= \int_T \E\left[|X''(t)|\mathbbm{1}_{\{X''(t) < 0\}} | X'(t)=0\right] p_{X'(t)}(0) dt\\
				&= \int_T \delta_t\left[\phi(b_t)+b_t\Phi(b_t)\right]\frac{1}{\sqrt{2\pi \la_1(t)}}e^{-\frac{m'(t)^2}{2\la_1(t)}} dt, 
			\end{split}
		\end{equation*} 
	yielding \eqref{eq:LocalMax-1D-1}. Next, we turn to computing $\E[M_u(X, T)]$. 
	
	Let $Y_1 = X(t)-m(t)$ and $Y_2=X''(t)/\delta_t + b_t$. It then follows from \eqref{eq:cov-conditional-m} that 
	\begin{equation}\label{eq:cov-conditional-m2}
		\begin{split}
			(	Y_1, Y_2 |  X'(t)=0) \sim N\left(\begin{pmatrix}
				0 \\
				0
			\end{pmatrix},\, \begin{pmatrix}
				1 & \rho_t\\
				\rho_t & 1
			\end{pmatrix}\right).
		\end{split}
	\end{equation} 
Then 
\begin{equation}\label{eq:G}
	\begin{split}
		G_t(u)&:= \E\left[|X''(t)|\mathbbm{1}_{\{X''(t) < 0\}}\mathbbm{1}_{\{X(t)>u\}} |  X'(t)=0\right]\\
		&= -\delta_t\E\left[\left(Y_2 -b_t \right)\mathbbm{1}_{\{Y_2 < b_t\}}\mathbbm{1}_{\{Y_1>u-m(t)\}} \Big|  X'(t)=0\right]\\
		&=  -\delta_t\int_{u-m(t)}^\infty dy_1 \int_{-\infty}^{b_t} \frac{y_2-b_t}{2\pi\sqrt{1-\rho_t^2}} e^{-\frac{y_1^2-2\rho_t y_1y_2 + y_2^2}{2(1-\rho_t^2)}}dy_2.
	\end{split}
\end{equation} 
Making change of variables $x=y_1/\sqrt{1-\rho_t^2}$ and $y=y_2/\sqrt{1-\rho_t^2}$, we obtain
\begin{equation*}
	\begin{split}
		G_t(u)& = -\frac{\delta_t(1-\rho_t^2)}{2\pi} \int_{\frac{u-m(t)}{\sqrt{1-\rho_t^2}}}^\infty dx \int_{-\infty}^{\frac{b_t}{\sqrt{1-\rho_t^2}}} \left(y- \frac{b_t}{\sqrt{1-\rho_t^2}}\right) e^{-\frac{x^2-2\rho_t xy + y^2}{2}}dy\\
		&= -\delta_t(1-\rho_t^2) \int_{\frac{u-m(t)}{\sqrt{1-\rho_t^2}}}^\infty \frac{1}{\sqrt{2\pi}} e^{-\frac{(1-\rho_t^2)x^2}{2}} dx \int_{-\infty}^{\frac{b_t}{\sqrt{1-\rho_t^2}}} \left(y- \frac{b_t}{\sqrt{1-\rho_t^2}}\right) \frac{1}{\sqrt{2\pi}} e^{-\frac{(y-\rho_t x)^2}{2}}dy\\
		&= -\delta_t(1-\rho_t^2)\int_{\frac{u-m(t)}{\sqrt{1-\rho_t^2}}}^\infty \frac{1}{\sqrt{2\pi}} e^{-\frac{(1-\rho_t^2)x^2}{2}} \Bigg[ \left(\rho_t x - \frac{b_t}{\sqrt{1-\rho_t^2}}\right) \Phi\left(\frac{b_t}{\sqrt{1-\rho_t^2}}-\rho_t x\right) \\
		&\qquad - \phi\left(\frac{b_t}{\sqrt{1-\rho_t^2}}-\rho_t x\right)   \Bigg] dx.
	\end{split}
\end{equation*} 
Applying integration by parts, we obtain
\begin{equation}\label{eq:1D-num-m2}
	\begin{split}
		G_t(u)& = \delta_t \Bigg[ \phi(b_t)\Psi\left(\frac{u-m(t)-\rho_tb_t}{\sqrt{1-\rho_t^2}}\right) - \rho_t\phi(u-m(t))\Phi\left(\frac{b_t-\rho_t(u-m(t))}{\sqrt{1-\rho_t^2}}\right)  \\
		&\qquad \quad  + b_t \int_{u-m(t)}^\infty \phi(z) \Phi\left(\frac{b_t - \rho_t z}{\sqrt{1-\rho_t^2}}\right) dz\Bigg].
	\end{split}
\end{equation} 
Note that the last integral can be written as 
\begin{equation}\label{eq:1D-num-m3}
	\begin{split}
		\int_{u-m(t)}^\infty \phi(z) \Phi\left(\frac{b_t - \rho_t z}{\sqrt{1-\rho_t^2}}\right) dz&= \P\left(\mathcal{N}_1>u-m(t), \mathcal{N}_2<\frac{b_t - \rho_t \mathcal{N}_1}{\sqrt{1-\rho_t^2}}\right) \\
		&= \P\left(\mathcal{N}_1>u-m(t), \rho_t \mathcal{N}_1 + \sqrt{1-\rho_t^2} \mathcal{N}_2 <b_t\right)\\
		&= \P(\xi_1>u-m(t), \xi_2 <b_t),
	\end{split}
\end{equation} 
where $\mathcal{N}_1$ and $\mathcal{N}_2$ are i.i.d. standard normal random variables, and $\xi_1$ and $\xi_2$ have a bivariate normal distribution as shown in \eqref{eq:xi}. Plugging \eqref{eq:1D-num-m2} and \eqref{eq:1D-num-m3} into the Kac-Rice formula 
\[
\E[M_u(X, T)] = \int_T G_t(u)\frac{1}{\sqrt{2\pi \la_1(t)}}e^{-\frac{m'(t)^2}{2\la_1(t)}} dt,
\]
we obtain \eqref{eq:LocalMax-1D-2}.
		\end{proof}

\begin{remark}
	As observed in Section \ref{sec:2D} below, similar to Gaussian fields in $\R^N$ with $N>1$, the usual method to compute $G_t(u)$ in \eqref{eq:G} is as follows:
	\begin{equation*}
		G_t(u) = \int_u^\infty \E\left[|X''(t)| \mathbbm{1}_{\{X''(t)< 0\}}|X(t)=x,  X'(t)=0\right] p_{X(t)}(x)dx.
	\end{equation*}
	However, this method of computation results in a complicated expression, making it difficult to extract the parameter $\rho_t$ for simplification. In particular, it would be difficult to determine that the peak height distribution for the centered case in \eqref{eq:PHD-1D-0} below depends solely on the parameter $\rho_t$.
\end{remark}

\begin{remark}
	For a centered Gaussian process, we have $m(t)\equiv 0$, implying $b_t\equiv 0$ and hence
	\begin{align}		
		\E[M(X, T)]&=\int_T \frac{\delta_t}{2\pi\sqrt{ \la_1(t)}} dt, \nonumber	\\
		\E[M_u(X, T)]&= \int_T \Bigg[ \frac{1}{\sqrt{2\pi}}\Psi\left(\frac{u}{\sqrt{1-\rho_t^2}}\right) - \rho_t\phi(u)\Phi\left(\frac{-\rho_tu}{\sqrt{1-\rho_t^2}}\right)   \Bigg]\frac{\delta_t}{\sqrt{2\pi \la_1(t)}} dt. \nonumber
	\end{align} 
Note also that, for the stationary case, since $\la_1'(t)\equiv 0$, the parameter $b_t$ simplifies to $b_t=-m''(t)/\delta_t=-m''(t)/\la_2(t)$.
\end{remark}

The results below present the derived peak height distributions. Notably, the peak height distribution of a centered nonstationary Gaussian process given in \eqref{eq:PHD-1D-0} depends only on the conditional correlation parameter $\rho_t$. 
\begin{corollary}\label{cor:1D_PHD}
	Let $\{X(t), t\in \R\}$ be a smooth, unit-variance, nonstationary Gaussian process with mean function $m(t)$. Then the peak height density at $t$ is given by 
	\begin{equation}\label{eq:PHD-1D}
		\begin{split}
			h_t(x)&= \left(\sqrt{1-\rho_t^2}\phi\left(\frac{b_t-\rho_t(x-m(t))}{\sqrt{1-\rho_t^2}}\right) + [b_t- \rho_t(x-m(t))]\Phi\left(\frac{b_t-\rho_t(x-m(t))}{\sqrt{1-\rho_t^2}}\right) \right)\\
			&\quad \times \phi(x-m(t)) [\phi(b_t)+b_t\Phi(b_t)]^{-1}, \quad x\in \R,
		\end{split}
	\end{equation} 
where $\rho_t\in (-1,0)$ and $b_t$ are as in \eqref{eq:rho} and \eqref{eq:xi} respectively. 
\end{corollary}
\begin{proof}
Since $X$ is a 1D Gaussian process, the peak height distribution in \eqref{eq:Palm distr Euclidean} becomes
\begin{equation}\label{eq:1D}
	\begin{split}
		F_t(u) 
		 = \frac{G_t(u)}{\delta_t\left[\phi(b_t)+b_t\Phi(b_t)\right]},
	\end{split}
\end{equation} 
where $G_t(u)$ is given in \eqref{eq:G} and the denominator due to \eqref{eq:J1-1D}. Taking the derivative with respect to $u$ in \eqref{eq:1D-num-m2}, we obtain that, $-G_t'(u)$ can be expressed as
\begin{equation*}
	\begin{split}
		&\quad \delta_t\Bigg[\frac{\phi(b_t)}{\sqrt{1-\rho_t^2}}\phi\left(\frac{u-m(t)-\rho_tb_t}{\sqrt{1-\rho_t^2}}\right) - \rho_t(u-m(t))\phi(u-m(t))\Phi\left(\frac{b_t-\rho_t[u-m(t)]}{\sqrt{1-\rho_t^2}}\right)  \\
		&\quad - \frac{\rho_t^2}{\sqrt{1-\rho_t^2}}\phi(u-m(t))\phi\left(\frac{b_t-\rho_t[u-m(t)]}{\sqrt{1-\rho_t^2}}\right) + b_t \phi(u-m(t)) \Phi\left(\frac{b_t-\rho_t[u-m(t)]}{\sqrt{1-\rho_t^2}}\right)\Bigg]\\
		&= \delta_t\phi(u-m(t))\Bigg[\sqrt{1-\rho_t^2}\phi\left(\frac{b_t-\rho_t[u-m(t)]}{\sqrt{1-\rho_t^2}}\right) + [b_t- \rho_t(u-m(t))]\Phi\left(\frac{b_t-\rho_t[u-m(t)]}{\sqrt{1-\rho_t^2}}\right)  \Bigg].
	\end{split}
\end{equation*} 
Plugging this into the derivative of \eqref{eq:1D} and noting the fact $h_t(u) = -F_t'(u)$, we obtain the peak height density in \eqref{eq:PHD-1D}. 
\end{proof}
	
	Here, we make some remarks on the comparison between the peak height distributions for certered nonstationary and stationary Gaussian processes in 1D. 
	\begin{remark}
		As an important result, note that, if $m(t)\equiv 0$, then $b_t\equiv 0$; plugging this into \eqref{eq:PHD-1D} yields the following peak height density for a centered nonstationary Gaussian process:
		\begin{equation}\label{eq:PHD-1D-0}
			\begin{split}
				h_t(x)&= \sqrt{1-\rho_t^2}\phi\left(\frac{x}{\sqrt{1-\rho_t^2}}\right) - \sqrt{2\pi}\rho_t x\phi(x)\Phi\left(\frac{-\rho_t x}{\sqrt{1-\rho_t^2}}\right), \quad x\in \R.
			\end{split}
		\end{equation} 
		To make comparison with the stationary case, let $X$ be a centered stationary Gaussian process. Then, the covariance depends only on the distance of two points, and hence can be written as ${\rm Cov}(X(t), X(s)) = \varphi\left(|t-s|^2\right)$, where $t,s\in \R$ and $\varphi: [0,\infty) \to \R$ is an appropriate function. The peak height density can be found in \cite{CL67} or derived from the general formula for Gaussian fields in \cite{CS15, CS18} as follows:
		\begin{equation}\label{eq:1D-stationary}
			\begin{split}
				h(x)=\sqrt{1-\kappa^2/3}\phi\left(\frac{ x}{\sqrt{1-\kappa^2/3}} \right) + \frac{ \sqrt{2\pi}\kappa}{\sqrt{3}} x\phi(x)\Phi\left(\frac{\kappa x}{\sqrt{3-\kappa^2}} \right), \quad x\in \R,
			\end{split}
		\end{equation}
		where $\kappa=-\varphi'(0)/\sqrt{\varphi''(0)}$, or $\kappa=\sqrt{3}{\rm Var}(X'(t))/\sqrt{{\rm Var}(X''(t))}$ since ${\rm Var}(X'(t))=-2\varphi'(0)$ and ${\rm Var}(X''(t))=12\varphi''(0)$ as shown in \cite{CS15}. Note that $\la_1'(t)=0$ due to the independence between $X'(t)$ and $X''(t)$ under stationarity, thus $\kappa=-\sqrt{3}\rho_t$ by \eqref{eq:rho}. Plugging this into \eqref{eq:1D-stationary}, we see that the two expressions for the peak height distribution, \eqref{eq:PHD-1D-0} and \eqref{eq:1D-stationary}, coincide for the stationary case. 
	\end{remark}
	
	\begin{remark}
		The peak height distribution in \eqref{eq:1D-stationary} has been used in computing p-values in statistical applications such as signal detection \cite{CS17} and change point detection \cite{ChangePoint}. The derived formulas in \eqref{eq:PHD-1D} and \eqref{eq:PHD-1D-0} will be valuable for conducting multiple tests for both signal and change point detection under nonstationary Gaussian noise. A notable advantage of this approach is that the peak height distribution in \eqref{eq:PHD-1D-0} relies solely on the conditional correlation parameter $\rho_t$, which simplifies parameter estimation and statistical inference. As demonstrated in \eqref{eq:rho}, estimating the parameter $\rho_t$ requires only the variances of the first and second derivatives of the Gaussian process.
	\end{remark}
	
	\begin{remark}
		By comparing \eqref{eq:1D-stationary} and \eqref{eq:PHD-1D-0}, we observe that, in the 1D case, the extension from stationarity to nonstationarity maintains the peak height distribution's dependence on a single parameter. It has been shown in \cite{CS15, CS18} that, for isotropic Gaussian fields in $\R^N$ with $N>1$, the peak height distribution also depends only on the parameter $\kappa$. However, this property typically does not hold beyond isotropy. As we will demonstrate in Section \ref{sec:2D}, even for 2D centered stationary nonisotropic Gaussian fields, the peak height distribution generally depends on at least two parameters.
	\end{remark}
	
	\subsection{The degenerate case: $\rho_t=-1$}
	In this section, we derive the results in Corollary \ref{cor:1D-1} below for the degenerate case $\rho_t=-1$. It is important to note that, due to the degeneracy, the arguments used to derive \eqref{eq:G} are no longer applicable. Instead, we will leverage the linear relation between $X(t)$ and $X''(t)$ conditional on $X'(t)=0$, which is implied by the degeneracy, to compute the conditional expectation in \eqref{eq:G}.  
	
	\begin{corollary}\label{cor:1D-1}
		Let $\{X(t), t\in \R\}$ be a smooth, unit-variance, nonstationary Gaussian process. Then the peak height density at $t$ with $\rho_t=-1$ is given by 
		\begin{equation}\label{eq:PHD-1D-1}
			\begin{split}
				h_t(x)&= \phi(x-m(t))[x-m(t)+b_t][\phi(b_t) + b_t\Phi(b_t) ]^{-1}\mathbbm{1}_{\{x> m(t)-b_t\}}.
			\end{split}
		\end{equation} 
	\end{corollary}
\begin{proof}
	Note that,  as a conditional correlation, $\rho_t=-1$ implies that there exist constants $c_1<0$ and $c_2\in \R$ such that 
	\begin{equation}\label{eq:1D-degen}
	(X(t) \, |\, X'(t)=0) = c_1(X''(t) \, |\, X'(t)=0) + c_2.
\end{equation} 
	Taking the mean and variance on both sides yields $m(t)=-c_1\delta_tb_t+c_2$ and $1=c_1^2\delta_t^2$, which implies $c_1=-1/\delta_t$ and $c_2=m(t)-b_t$. Plugging \eqref{eq:1D-degen} into $G_t(u)$ in \eqref{eq:G}, we obtain that, for $u>m(t)-b_t$,
	\begin{equation*}
		\begin{split}
			G_t(u) &= \E\left[|X''(t)|\mathbbm{1}_{\{X''(t) < 0\}}\mathbbm{1}_{\{X(t)>u\}} |  X'(t)=0\right] \\
			&= \E\left[|X''(t)|\mathbbm{1}_{\{X''(t) < 0\}}\mathbbm{1}_{\{-X''(t)/\delta_t +m(t)-b_t>u\}} |  X'(t)=0\right]\\
			&= \E\left[|X''(t)|\mathbbm{1}_{\{X''(t)<\delta_t [m(t)-b_t-u]\}} |  X'(t)=0\right]\\
			&= -\int_{-\infty}^{\delta_t [m(t)-b_t-u]}\frac{x}{\sqrt{2\pi}\delta_t}e^{-\frac{(x+\delta_tb_t)^2}{2\delta_t^2}}dt\\
			&= \delta_t[\phi(u-m(t)) + b_t\Psi(u-m(t)) ];
		\end{split}
	\end{equation*} 
while for $u\le m(t)-b_t$, 
\[
G_t(u) = \E\left[|X''(t)|\mathbbm{1}_{\{X''(t) < 0\}} |  X'(t)=0\right] = \delta_t[\phi(b_t) + b_t\Phi(b_t) ],
\]
which is independent of $u$. Then, similarly to the proof of Corollary \ref{cor:1D_PHD}, taking the derivative $G_t'(u)$ and applying \eqref{eq:1D}, we obtain the desired result \eqref{eq:PHD-1D-1}.
\end{proof}

Note that, for the degenerate case $\rho_t=-1$, the support of the peak height density is $x>m(t)-b_t$,  which differs from the entire real line $\R$ for the nondegenerate case $-1<\rho_t<0$. It is also interesting to observe that the derived degenerate density in \eqref{eq:PHD-1D-1} coincides with the nondegenerate density in \eqref{eq:PHD-1D} when taking the limit $\rho_t\downarrow -1$.

A simple example of the nondegenerate case is the cosine process in $\R$. Specifically, let $X(t)=Z(t)+m(t)$ with
  \begin{equation}\label{eq:cosine}
	Z(t) = \mathcal{N}_1 \cos(\omega t) + \mathcal{N}_2 \sin(\omega t), \quad t\in \R,
\end{equation}
where $\mathcal{N}_1$ and $\mathcal{N}_2$ are independent standard normal random variables, and $\omega$ is a positive constant. Then $Z$ is a centered, stationary, unit-variance, smooth Gaussian process with covariance $\E[Z(t)Z(0)] = \cos(\omega t)$. The nondegeneracy arises from the fact that $Z''(t)=-\omega^2 Z(t)$, yielding $\rho_t\equiv -1$.

We can also observe from the proof of Corollary \ref{cor:1D-1} that, $\E[M(X, T)]$ remains the same as in \eqref{eq:LocalMax-1D-1} for all $\rho_t\in[-1,0)$; and if $\rho_t\equiv -1$, then
\begin{equation*}
	\begin{split}
		\E[M_u(X, T)]&= \begin{cases}
			\int_T \left[\phi(u-m(t)) + b_t\Psi(u-m(t)) \right]\frac{\delta_t}{\sqrt{2\pi \la_1(t)}} e^{-\frac{m'(t)^2}{2\la_1(t)}}dt, \quad &u>m(t)-b_t,\\
			\int_T \left[\phi(b_t)+b_t\Phi(b_t)\right]\frac{\delta_t}{\sqrt{2\pi \la_1(t)}} e^{-\frac{m'(t)^2}{2\la_1(t)}}dt, \quad  &u\le m(t)-b_t.
		\end{cases}
	\end{split}
\end{equation*}

	\section{Planar Stationary Gaussian Random Fields}\label{sec:2D}
	In this section, we study the expected number and height distribution of local maxima for stationary Gaussian fields on $\R^2$. Generally, these are difficult to evaluate due to the complexity of the Hessian matrix. In particular, without the isotropic property, the random matrix techniques based on GOE and GOI \cite{CS15,CS18} are not applicable. To obtain exact formulas, we impose an independence assumption on partial derivatives along distinct directions (see \eqref{eq:2D-condition} below). We leverage the fact that $\nabla^2 X(t)\prec 0$ is equivalent to both diagonal entries being negative and the determinant being positive. The main technique involves transforming the distribution of the determinant of $\nabla^2 X(t)$ into a generalized chi-squared distribution, allowing us to handle the expectation involving $\nabla^2 X(t)$. The exact formulas, derived in Theorem \ref{theorem:2D}, are expressed in terms of the generalized chi-squared density.
	
	\subsection{Notations and assumptions}
	For a function $f(\cdot) \in C^2(\R^2)$ and $t\in \R^2$, we introduce the following notations for derivatives:
	\begin{equation}\label{Eq:notatoin-diff}
		\begin{split}
			f_i (t)&=\frac{\partial f(t)}{\partial t_i}, \quad f_{jk}(t)=\frac{\partial^2 f(t)}{\partial t_j\partial t_k}, \quad \forall i, j, k=1, 2;\\
			 \nabla f(t) &=(f_1(t), f_2(t)), \quad \nabla^2 f(t) = \left(f_{\ell n}(t)\right)_{ \ell, n = 1, 2}.
		\end{split}
	\end{equation}
	Let $\{X(t): t\in T \subset \R^2\}$ be a smooth, unit-variance, stationary planar Gaussian random field with mean function $m(t)$. Specifically, following Corollary 11.3.2 in \cite{AT07}, we assume that the distributions of $(\nabla X(t), X_{11}(t), X_{22}(t), X_{12}(t))$ are nondegenerate for all $t\in T$, and that there exist constants $K, \alpha>0$ such that for all $t,s\in T$,
	\[
	\max_{i,j\in \{1, 2\}}|{\rm Var}(X_{ij}(t))-{\rm Cov}(X_{ij}(t), X_{ij}(s))| \le K|\log \|t-s\||^{-(1+\alpha)}.
	\]
	 We can write $X(t) = Z(t) + m(t)$, where $Z(t)$ is a centered stationary Gaussian field. Note that the distribution of $(Z(t), Z_i(t), Z_{jk}(t))$ does not depend on $t$ due to the stationarity of $Z$. We have the following spectral representations (see \cite[p.~112]{AT07}) for the covariance of the gradient and Hessian,
	\begin{equation}\label{eq:spectral-rep}
		\E[Z_i(t)Z_j(t)] = \int_{\R^2} x_ix_j \nu(dx), \quad \E[Z_{ij}(t)Z_{kl}(t)] = \int_{\R^4} x_ix_jx_kx_l \nu(dx), \quad i,j,k,l=1,2,
	\end{equation}
	where $\nu(\cdot)$ is a spectral measure. To obtain the exact formula, we will make use of the following independence assumption on partial derivatives along distinct directions:
	\begin{equation}\label{eq:2D-condition}
		\E[Z_1(t)Z_2(t)] = \E[Z_{11}(t)Z_{12}(t)]=\E[Z_{22}(t)Z_{12}(t)]=0.
	\end{equation} 
	Here is an example of a covariance satisfying this property. Let $C(t)=\E[Z(t)Z(0)]$, where $t=(t_1, t_2)\in \R^2$, such that
	\[
	C(t_1, t_2) = C(-t_1, t_2) \quad \text{ and } \quad C(t_1, t_2) = C(t_1,-t_2).
	\] 
	This indicates that, the spectral measure $\nu$ in \eqref{eq:spectral-rep} is symmetric in each direction $\lambda_1$ and $\lambda_2$, separately, a property known as quadrant symmetry. Such a quadrant symmetry implies \eqref{eq:2D-condition}.
	
	Note that, by \eqref{eq:spectral-rep}, ${\rm Cov}(X_{11}(t), X_{22}(t))={\rm Var}(X_{12}(t))$. Hence, under our assumption \eqref{eq:2D-condition}, the variance-covariance matrix of the Hessian can be written as
	\begin{equation}\label{eq:cov-Hessian}
		\begin{split}
			{\rm Cov}(X_{11}(t), X_{22}(t), X_{12}(t)) = \begin{pmatrix} 
				\sigma_{11}^2 & \sigma_{12}^2 & 0 \\
				\sigma_{12}^2 & \sigma_{22}^2 & 0 \\
				0 & 0 & \sigma_{12}^2 \\
			\end{pmatrix}.
		\end{split}
	\end{equation}
	
	\subsection{Expected number and height distribution of local maxima}
	First, we derive the following result, which will be useful for simplifying the conditions. 
	\begin{lemma}\label{lemma:2D}
		Let $\{\widetilde{X}(t), t=(t_1,t_2)\in \R^2\}$ be a smooth, unit-variance, stationary Gaussian field satisfying ${\rm Cov}(\nabla \widetilde{X}(t)) = {\rm diag}(\gamma_1^2, \gamma_2^2)$ and ${\rm Cov}(\widetilde{X}_{ii}(t), \widetilde{X}_{ij}(t))=0$, $i\neq j$. Let $X(t_1, t_2)= \widetilde{X}(t_1/\gamma_1, t_2/\gamma_2)$. Then $X$ is stationary with ${\rm Cov}(\nabla X(t)) = I_2$ and ${\rm Cov}(X_{ii}(t), X_{ij}(t))=0$, $i\neq j$. Moreover, $\E[M_u(X, T)] =  \E[M_u(\widetilde{X}, \widetilde{T})]$ and $\E[M(X, T)] = \E[M(\widetilde{X}, \widetilde{T})]$, where $\widetilde{T} = \{(t_1/\gamma_1, t_2/\gamma_2): (t_1, t_2)\in T\}$; and the peak height distribution of $X$ at $(t_1, t_2)$ is the same as that of $\widetilde{X}$ at $(t_1/\gamma_1, t_1/\gamma_2)$.
	\end{lemma}
	\begin{proof}
		Note that, the transformation is linear, and thus keeps the stationarity of $X$. Since $X(t_1, t_2)=\widetilde{X}(t_1/\gamma_1, t_2/\gamma_2)$, one has 
		\begin{equation*}
			\begin{split}
				X_i(t_1, t_2) &= \widetilde{X}_i(t_1/\gamma_1, t_2/\gamma_2)/\gamma_i, \quad i=1, 2,\\
				X_{jk}(t_1, t_2) &= \widetilde{X}_{jk}(t_1/\gamma_1, t_2/\gamma_2)/(\gamma_j\gamma_k), \quad j,k=1,2.
			\end{split}
		\end{equation*}
		Therefore, ${\rm Cov}(\nabla X(t)) = I_2$ and ${\rm Cov}(X_{ii}(t), X_{ij}(t))=0$, $i\neq j$. Moreover, we have 
		\begin{equation*}
			\begin{split}
				M_u(X, T) &= \# \{ (t_1, t_2)\in T: X(t_1, t_2)\geq u, \nabla X(t_1, t_2)=0, \nabla^2 X(t_1, t_2)\prec 0 \}\\
				&= \# \{ (t_1, t_2)\in T: \widetilde{X}(t_1/\gamma_1, t_2/\gamma_2)\geq u, \nabla \widetilde{X}(t_1/\gamma_1, t_2/\gamma_2)=0, \nabla^2 \widetilde{X}(t_1/\gamma_1, t_2/\gamma_2)\prec 0 \}\\
				&= M_u(\widetilde{X}, \widetilde{T}),
			\end{split}
		\end{equation*}
		yielding $\E[M_u(X, T)] =  \E[M_u(\widetilde{X}, \widetilde{T})]$. Similarly, we have $\E[M(X, T)] =  \E[M(\widetilde{X}, \widetilde{T})]$ and that the peak height distribution of $X$ at $(t_1, t_2)$ is the same as that of $\widetilde{X}$ at $(t_1/\gamma_1, t_1/\gamma_2)$.
	\end{proof}
	
	It is seen from Lemma \ref{lemma:2D} that, for a stationary planar Gaussian field $\{\widetilde{X}(t_1, t_2), (t_1, t_2)\in \widetilde{T}\}$ satisfying \eqref{eq:2D-condition}, the expected number and height distribution of local maxima can be computed through $\{X(t_1, t_2), (t_1, t_2)\in T\}$ which also satisfies \eqref{eq:2D-condition} with ${\rm Cov}(\nabla X(t)) = I_2$. Therefore, without loss of generality, we assume the covariance of the gradient of the Gaussian field is a $2\times 2$ identity matrix. 
	
	As mentioned earlier, we will need to use the densities of generalized chi-squared distributions. However, such densities typically do not have simple closed-form expressions. The following lemma provides a power series expansion for these generalized chi-squared densities.
	
	\begin{lemma}[Kotz et al. \cite{Kotz67}]\label{lemma:chi-square}
		Let $W=\alpha_1 (\mathcal{N}_1+\beta)^2 + \alpha_2 \mathcal{N}_2^2$, where $\mathcal{N}_1$ and $\mathcal{N}_2$ are independent standard normal random variables. Then the probability density function of $W$ is given by
		\begin{equation*}
			\begin{split}
				f_W(w) = \sum_{k=0}^\infty \frac{c_k(\alpha_1, \alpha_2, \beta)(-w)^k}{2^{k+1} k!},
			\end{split}
		\end{equation*}
		where $c_k(\alpha_1, \alpha_2, \beta)$ are determined by the expansion
		\begin{equation}\label{eq:ck}
			\begin{split}
				\sum_{k=0}^\infty c_k(\alpha_1, \alpha_2, \beta)\theta^k = (\alpha_1\alpha_2)^{-1/2}\exp\left\{-\frac{\beta^2}{2(1-\theta/\alpha_1)}\right\}[(1-\theta/\alpha_1)(1-\theta/\alpha_2)]^{-1/2}.
			\end{split}
		\end{equation}
	\end{lemma}
\begin{remark}\label{remark:chi}
	Note that, if $\alpha_1=\alpha_2=1$ and $\beta=0$, then $c_k(1,1,0)=1$ for all $k\ge 0$, so that \begin{equation*}
		\begin{split}
			f_W(w) = \frac{1}{2}\sum_{k=0}^\infty \frac{(-w)^k}{2^kk!}=\frac{1}{2}e^{-\frac{w}{2}},
		\end{split}
	\end{equation*}
	which is exactly the pdf of $\chi_2^2$, chi-squared distribution with degree of freedom 2.
	
	For the case of $\beta=0$, note that 
	\begin{equation*}
		\begin{split}
			(1-\theta/\alpha_1)^{-1/2} &= \sum_{i=0}^\infty {-1/2 \choose i} \left(-\frac{\theta}{\alpha_1}\right)^i = \frac{1}{\sqrt{\pi}}\sum_{i=0}^\infty \frac{\Gamma\left(i+\frac{1}{2}\right)}{i!}\frac{\theta^i}{\alpha_1^i},\\
			(1-\theta/\alpha_2)^{-1/2} &= \sum_{i=0}^\infty {-1/2 \choose i} \left(-\frac{\theta}{\alpha_2}\right)^i = \frac{1}{\sqrt{\pi}}\sum_{i=0}^\infty \frac{\Gamma\left(i+\frac{1}{2}\right)}{i!}\frac{\theta^i}{\alpha_2^i}.
		\end{split}
	\end{equation*}
	So, the coefficients for the Taylor series of $(\alpha_1\alpha_2)^{-1/2}[(1-\theta/\alpha_1)(1-\theta/\alpha_2)]^{-1/2}$ are given by
	\[
	c_k(\alpha_1, \alpha_2, 0) = \frac{1}{\pi\sqrt{\alpha_1\alpha_2}}\sum_{i=0}^k \frac{\Gamma\left(i+\frac{1}{2}\right)\Gamma\left(k-i+\frac{1}{2}\right)}{i!(k-i)!\alpha_1^i\alpha_2^{k-i}}.
	\]
\end{remark}

	The following result is a direct consequence of integration by parts.
	\begin{lemma} Let $j$ be a nonnegative integer. Then 
		\begin{equation}\label{eq:Q}
			\begin{split}
		Q_j(y):=\int_{-\infty}^y s^j \phi(s)ds= \begin{cases}
			\Phi(y) & \text{\rm if } j=0,\\
			-\sum_{i=0}^n \frac{(2n)!!}{(2i)!!}y^{2i}\phi(y) & \text{\rm if } j=2n + 1,\\
			(2n+1)!!\Phi(y) - \sum_{i=0}^n \frac{(2n+1)!!}{(2i+1)!!}y^{2i+1}\phi(y)  & \text{\rm if } j=2n+2,
		\end{cases}
	\end{split}
\end{equation}
		where $n=0, 1, 2, \ldots$
	\end{lemma}

To address the derived results in a more concise way, we introduce the following notations:
\begin{equation}\label{eq:a_t}
	\begin{split}
		\rho &= \frac{\sigma_{12}^2}{\sigma_{11}\sigma_{22}}, \quad  \tilde{\rho} = \frac{\tilde{\sigma}_{12}^2}{\tilde{\sigma}_{11}\tilde{\sigma}_{22}}, \quad \tilde{\sigma}_{ij}^2 = \sigma_{ij}^2-1, \quad i,j=1,2, \\
		a_t &= \frac{m_{11}(t)}{\sigma_{11}}+\frac{m_{22}(t)}{\sigma_{22}}, \quad \tilde{a}_t(x)=\frac{m_{11}(t)-x}{\tilde{\sigma}_{11}}+\frac{m_{22}(t)-x}{\tilde{\sigma}_{22}},\\
		c_{k,t}&=c_k\left( \frac{1-\rho}{2},\, \rho, \,  \frac{1}{\sqrt{2(1-\rho)}}\left(\frac{m_{11}(t)}{\sigma_{11}}-\frac{m_{22}(t)}{\sigma_{22}}\right) \right),\\
		\tilde{c}_{k,t}(x)&= c_k\left( \frac{1-\tilde{\rho}}{2},\,  \tilde{\rho}, \,  \frac{1}{\sqrt{2(1-\tilde{\rho})}}\left(\frac{m_{11}(t)-x}{\tilde{\sigma}_{11}}-\frac{m_{22}(t)-x}{\tilde{\sigma}_{22}}\right) \right),
	\end{split}
\end{equation}
where $\sigma_{11}^2$, $\sigma_{22}^2$ and $\sigma_{12}^2$ are given in \eqref{eq:cov-Hessian}, and the constants $c_k(\cdot, \cdot, \cdot)$ are given in Lemma \ref{lemma:chi-square}. The following main results provide the exact formulas for the expected number and height distributions of local maxima for planar stationary Gaussian fields with a mean function.
	
	\begin{theorem}\label{theorem:2D}
		Let $\{X(t), t\in T \subset \R^2\}$ be a smooth, unit-variance, stationary Gaussian field with mean function $m(t)$ and satisfying ${\rm Cov}(\nabla X(t)) = I_2$, ${\rm Cov}(X_{ii}(t), X_{ij}(t))=0$, $i\neq j$, and $\nabla^2 m(t)={\rm diag}(m_{11}(t), m_{22}(t))$. Then the expected number and height distribution of local maxima are given respectively by
		\begin{align}		
			\E[M(X, T)]&=\frac{1}{2\pi}\int_T e^{-\frac{m_1^2(t)+m_2^2(t)}{2}}J_{1,t} \, dt, \label{eq:LocalMax-2D-1}\\
			\E[M_u(X, T)]&= \frac{1}{(2\pi)^{3/2}}\int_T\int_u^\infty  e^{-\frac{(x-m(t))^2+m_1^2(t)+m_2^2(t)}{2}} J_{2,t}(x)\, dxdt, \label{eq:LocalMax-2D-2}\\
			h_t(x)&= \phi(x-m(t))J_{2,t}(x)/J_{1,t}, \label{eq:phd1},
		\end{align} 
		where 
		\begin{align}		
			J_{1,t}&=\sigma_{11}\sigma_{22}\sum_{k=0}^\infty \frac{(-1)^{k+1}a_t^{2k+4}c_{k,t}}{2^{3k+5}(k+2)!} \sum_{j=0}^{2k+4} {2k+4 \choose j} \left(\frac{\sqrt{2(1+\rho)}}{a_t}\right)^jQ_j\left(\frac{-a_t}{\sqrt{2(1+\rho)}}\right), \label{eq:J1t}	\\
			J_{2,t}(x)&= \tilde{\sigma}_{11}\tilde{\sigma}_{22}\sum_{k=0}^\infty \frac{(-1)^{k+1}\tilde{a}_t^{2k+4}(x)\tilde{c}_{k,t}(x)}{2^{3k+5}(k+2)!} \sum_{j=0}^{2k+4} {2k+4 \choose j}  \left(\frac{\sqrt{2(1+\tilde{\rho})}}{\tilde{a}_t(x)}\right)^jQ_j\left(\frac{-\tilde{a}_t(x)}{\sqrt{2(1+\tilde{\rho})}}\right), \label{eq:J2t}
		\end{align} 
		and $\sigma_{ij}$, $Q_j(\cdot)$, $\tilde{\sigma}_{ij}$, $\rho$, $\tilde{\rho}$, $a_t$, $\tilde{a}_t(x)$, $c_{k,t}$ and $\tilde{c}_{k,t}(x)$ are given in \eqref{eq:cov-Hessian}, \eqref{eq:Q} and \eqref{eq:a_t}.
	\end{theorem}
	\begin{proof}
		Note that, by \eqref{eq:cov-Hessian},  we have ${\rm Var}(X_{ij}(t))=\sigma_{ij}^2$, $i,j=1,2$. Let
		\begin{equation}\label{eq:J-def}
			\begin{split}
				J_{1,t} &= \E\left[|{\rm det} (\nabla^2 X(t))| \mathbbm{1}_{\{\nabla^2 X(t)\prec 0\}}\right],\\ 
		J_{2,t}(x) &= \E\left[|{\rm det} (\nabla^2 X(t))| \mathbbm{1}_{\{\nabla^2 X(t)\prec 0\}} | X(t)=x\right].
	\end{split}
\end{equation}
		Note that, a symmetric $2\times 2$ matrix is negative definite if and only if both diagonal entries are negative and the determinant is positive. Using this fact, we can write 
		\begin{equation}\label{eq:J1}
			\begin{split}
				J_{1,t} = \E \left[ \left(X_{11}(t)X_{22}(t)-X_{12}^2(t)\right)\mathbbm{1}_{\{X_{11}(t)<0\}}\mathbbm{1}_{\{X_{22}(t)<0\}}\mathbbm{1}_{\{X_{11}(t)X_{22}(t)-X_{12}^2(t)>0\}} \right].
			\end{split}
		\end{equation}
		Recall $X(t)=Z(t) + m(t)$. Let $V_1=Z_{11}(t)/\sigma_{11}$, $V_2=Z_{22}(t)/\sigma_{22}$ and $V_3=Z_{12}(t)/\sigma_{12}$. Then
		\begin{equation*}
			\begin{split}
				(V_1, V_2, V_3) 
				\sim N\left(\begin{pmatrix} 
					0 \\
					0 \\
					0 
				\end{pmatrix},  \begin{pmatrix} 
					1 & \rho & 0 \\
					\rho & 1 & 0 \\
					0 & 0 & 1 
				\end{pmatrix} \right).
			\end{split}
		\end{equation*}
		Let $Y_1 = V_1$, $Y_2 = V_3$ and $Y_3 = V_1 + V_2$. Then $Y_3 \sim N(0, 2(1+\rho))$ and 
		\begin{equation}\label{eq:Y1-Y2}
			\begin{split}
				(Y_1, Y_2 | Y_3 = y) \sim N\left(\begin{pmatrix} 
					\frac{y}{2} \\
					0 
				\end{pmatrix},  \begin{pmatrix} 
					\frac{1-\rho}{2} & 0 \\
					0 & 1 
				\end{pmatrix} \right).
			\end{split}
		\end{equation}
	Moreover, 
	\begin{equation}\label{eq:X_ij-tran}
		\begin{split}
			X_{11}(t) = \sigma_{11}Y_1+m_{11}(t), \ X_{22}(t)=\sigma_{22}(Y_3-Y_1)+m_{22}(t), \ X_{12}(t) = \sigma_{12}Y_2.
		\end{split}
	\end{equation}
		By \eqref{eq:Y1-Y2} and \eqref{eq:X_ij-tran}, we can write the conditional distribution of $(X_{11}(t)X_{22}(t)-X_{12}^2(t) | Y_3 = y)$ as
		\begin{equation}\label{eq:det-conditional}
			\begin{split}
				&\quad \left( (\sigma_{11}Y_1+m_{11}(t))[\sigma_{22}(Y_3-Y_1)+m_{22}(t)] - \sigma_{12}^2Y_2^2 \,|\, Y_3 = y \right)\\
				&\overset{d}{=}\left[\sigma_{11}\left(\sqrt{\frac{1-\rho}{2}}\mathcal{N}_1+\frac{y}{2}\right)+m_{11}(t)\right]\left[-\sigma_{22}\left(\sqrt{\frac{1-\rho}{2}}\mathcal{N}_1-\frac{y}{2}\right)+m_{22}(t)\right] - \sigma_{12}^2\mathcal{N}_2^2\\
				&\overset{d}{=} -\sigma_{11}\sigma_{22}\left[\frac{1-\rho}{2}\mathcal{N}_1^2 +\sqrt{\frac{1-\rho}{2}}\left(\frac{m_{11}(t)}{\sigma_{11}}-\frac{m_{22}(t)}{\sigma_{22}}\right)\mathcal{N}_1\right] \\
				&\quad + \left(\frac{\sigma_{11}y}{2}+m_{11}(t)\right)\left(\frac{\sigma_{22}y}{2}+m_{22}(t)\right)- \sigma_{12}^2\mathcal{N}_2^2\\
				&\overset{d}{=} -\sigma_{11}\sigma_{22}\Bigg[\frac{1-\rho}{2} \left(\mathcal{N}_1 + \frac{1}{\sqrt{2(1-\rho)}}\left(\frac{m_{11}(t)}{\sigma_{11}}-\frac{m_{22}(t)}{\sigma_{22}}\right)\right)^2 +\rho \mathcal{N}_2^2- \frac{1}{4}\left(y+a_t\right)^2 \Bigg],
			\end{split}
		\end{equation}
		where $a_t$ is from \eqref{eq:a_t}, and $\mathcal{N}_1$ and $\mathcal{N}_2$ are independent standard normal variables. Let 
		\begin{equation}\label{eq:W}
		W = \frac{1-\rho}{2} \left(\mathcal{N}_1 + \frac{1}{\sqrt{2(1-\rho)}}\left(\frac{m_{11}(t)}{\sigma_{11}}-\frac{m_{22}(t)}{\sigma_{22}}\right)\right)^2 +\rho \mathcal{N}_2^2.
	\end{equation}
		Then $W$ has a generalized chi-squared distribution. Note that 
		\begin{equation*}
			\begin{split}
		&\quad X_{11}(t)<0, \, X_{22}(t)<0 \, {\rm\text{ and }} \, X_{11}(t)X_{22}(t)-X_{12}^2(t)>0 \\ &\Leftrightarrow X_{11}(t)/\sigma_{11} + X_{22}(t)/\sigma_{22}<0\, {\rm\text{ and }} \, X_{11}(t)X_{22}(t)-X_{12}^2(t)>0.
	\end{split}
\end{equation*}
		Plugging this into \eqref{eq:J1}, together with \eqref{eq:X_ij-tran}, \eqref{eq:det-conditional} and the fact $Y_3 \sim N(0, 2(1+\rho))$, we obtain
		\begin{equation}\label{eq:J1E}
			\begin{split}
				J_{1,t} &= \E \Bigg[ \left((\sigma_{11}Y_1+m_{11}(t))[\sigma_{22}(Y_3-Y_1)+m_{22}(t)] - \sigma_{12}^2Y_2^2\right)\mathbbm{1}_{\{Y_3<-\frac{m_{11}(t)}{\sigma_{11}}-\frac{m_{22}(t)}{\sigma_{22}}\}}\\
				&\qquad \ \times \mathbbm{1}_{\{(\sigma_{11}Y_1+m_{11}(t))[\sigma_{22}(Y_3-Y_1)+m_{22}(t)] - \sigma_{12}^2Y_2^2>0\}} \Bigg]\\
				&= \frac{-\sigma_{11}\sigma_{22}}{2\sqrt{\pi(1+\rho)}} \int_{-\infty}^{-a_t} \E \left[ \left(W- \frac{1}{4}\left(y+a_t\right)^2\right) \mathbbm{1}_{\{W< \frac{1}{4}\left(y+a_t\right)^2\}}\right] e^{-\frac{y^2}{4(1+\rho)}}dy.
			\end{split}
		\end{equation}
	Applying Lemma \ref{lemma:chi-square} with $\alpha_1=\frac{1-\rho}{2}$, $\alpha_2=\rho$ and $\beta = \frac{1}{\sqrt{2(1-\rho)}}\left(\frac{m_{11}(t)}{\sigma_{11}}-\frac{m_{22}(t)}{\sigma_{22}}\right)$, we have the pdf of $W$ as follows:
	\begin{equation*}
		\begin{split}
			f_{W}(w) = \sum_{k=0}^\infty \frac{(-1)^kc_{k,t}w^k}{2^{k+1}k!},
		\end{split}
	\end{equation*}
	where $c_{k,t}$ are given in \eqref{eq:a_t}. This implies
	\begin{equation*}
		\begin{split}
			\E \left[ W\mathbbm{1}_{\{W< \frac{1}{4}\left(y+a_t\right)^2\}} \right] &= \sum_{k=0}^\infty \frac{(-1)^kc_{k,t}\left(y+a_t\right)^{2k+4}}{2^{k+1}k!4^{k+2}(k+2)},\\
			\E \left[ \mathbbm{1}_{\{W< \frac{1}{4}\left(y+a_t\right)^2\}} \right] &= \sum_{k=0}^\infty \frac{(-1)^kc_{k,t}\left(y+a_t\right)^{2k+2}}{2^{k+1}(k+1)!4^{k+1}},
		\end{split}
	\end{equation*}
	and hence 
	\begin{equation*}
		\begin{split}
			&\quad \E \left[ \left(W- \frac{1}{4}\left(y+a_t\right)^2\right) \mathbbm{1}_{\{W< \frac{1}{4}\left(y+a_t\right)^2\}}\right]= \sum_{k=0}^\infty \frac{(-1)^kc_{k,t}\left(y+a_t\right)^{2k+4}}{2^{k+1}4^{k+2}(k+2)!}\\
			&=\sum_{k=0}^\infty \frac{(-1)^kc_{k,t}}{2^{3k+5}(k+2)!}\sum_{j=0}^{2k+4} {2k+4 \choose j}a_t^{2k+4-j}y^j.
		\end{split}
	\end{equation*}
	Making change of variables $s=y/\sqrt{2(1+\rho)}$, we obtain
	\begin{equation*}
		\begin{split}
			J_{1,t}&=\frac{-\sigma_{11}\sigma_{22}}{2\sqrt{\pi(1+\rho)}} \sum_{k=0}^\infty \frac{(-1)^kc_{k,t}}{2^{3k+5}(k+2)!}\sum_{j=0}^{2k+4} {2k+4 \choose j}a_t^{2k+4-j}\int_{-\infty}^{-a_t} y^j e^{-\frac{y^2}{4(1+\rho)}}dy\\
			&=\frac{-\sigma_{11}\sigma_{22}}{2\sqrt{\pi(1+\rho)}} \sum_{k=0}^\infty \frac{(-1)^kc_{k,t}}{2^{3k+5}(k+2)!}\sum_{j=0}^{2k+4} {2k+4 \choose j}a_t^{2k+4-j}(\sqrt{2(1+\rho)})^{j+1}\int_{-\infty}^{\frac{-a_t}{\sqrt{2(1+\rho)}}} s^j e^{-\frac{s^2}{2}}ds\\
			&=\sigma_{11}\sigma_{22} \sum_{k=0}^\infty \frac{(-1)^{k+1}a_t^{2k+4}c_{k,t}}{2^{3k+5}(k+2)!}\sum_{j=0}^{2k+4} {2k+4 \choose j}\left(\frac{\sqrt{2(1+\rho)}}{a_t}\right)^jQ_j\left(\frac{-a_t}{\sqrt{2(1+\rho)}}\right),
		\end{split}
	\end{equation*}
		where $Q_j(\cdot)$ is given by \eqref{eq:Q}. Note that $p_{\nabla X(t)}(0,0) = \frac{1}{2\pi}e^{-\frac{m_1^2(t)+m_2^2(t)}{2}}$, and by stationarity,
		\[
		\E\left[|{\rm det} (\nabla^2 X(t))| \mathbbm{1}_{\{\nabla^2 X(t)\prec 0\}} \big| \nabla X(t)=(0, 0) \right] = \E\left[|{\rm det} (\nabla^2 X(t))| \mathbbm{1}_{\{\nabla^2 X(t)\prec 0\}}\right]=J_{1,t}.
		\]
		Then, \eqref{eq:LocalMax-2D-1} follows immediately from the Kac-Rice formula 
		\begin{equation*}
			\begin{split}
		\E[M(X,T)]&=\int_T \E\left[|{\rm det} (\nabla^2 X(t))| \mathbbm{1}_{\{\nabla^2 X(t)\prec 0\}} \big| \nabla X(t)=(0, 0) \right]p_{\nabla X(t)}(0,0) dt.
	\end{split}
\end{equation*}

Next, we turn to computing $J_{2,t}(x)$. Note that
\begin{equation*}
	\left(Z_{11}(t),Z_{22}(t),Z_{12}(t)|X(t) = x\right) \sim N \left( \begin{pmatrix}
		-x\\
		-x\\
		0
	\end{pmatrix},
	\begin{pmatrix}
		\sigma_{11}^2 - 1&\sigma_{12}^2-1&0\\
		\sigma_{12}^2 -1&\sigma_{22}^2-1&0\\
		0&0&\sigma_{12}^2
	\end{pmatrix}
	\right).
\end{equation*}
Recall $\tilde{\sigma}_{ij}^2 = \sigma_{ij}^2-1$.
Let $\tilde{V}_1=(Z_{11}(t)+x)/\tilde{\sigma}_{11}$, $\tilde{V}_2=(Z_{22}(t)+x)/\tilde{\sigma}_{22}$ and $\tilde{V}_3=Z_{12}(t)/\sigma_{12}$. Then
\begin{equation*}
	\begin{split}
		(\tilde{V}_1, \tilde{V}_2, \tilde{V}_3 | X(t)=x)
		\sim N\left(\begin{pmatrix} 
			0 \\
			0 \\
			0 \\
		\end{pmatrix},  \begin{pmatrix} 
			1 & \tilde{\rho} & 0 \\
			\tilde{\rho} & 1 & 0 \\
			0 & 0 & 1 \\
		\end{pmatrix} \right),
	\end{split}
\end{equation*}
where $\tilde{\rho}$ is given in \eqref{eq:a_t}. Let $\tilde{Y}_1 = \tilde{V}_1$, $\tilde{Y}_2 = \tilde{V}_3$ and $\tilde{Y}_3 = \tilde{V}_1 + \tilde{V}_2$. Then $(\tilde{Y}_3|X(t)=x) \sim N(0, 2(1+\tilde{\rho}))$ and
\begin{equation*}
	\begin{split}
		&(\tilde{Y}_1, \tilde{Y}_2 | X(t)=x, \tilde{Y}_3 = y) \sim N\left(\begin{pmatrix} 
			\frac{y}{2} \\
			0 
		\end{pmatrix},  \begin{pmatrix} 
			\frac{1-\tilde{\rho}}{2} & 0 \\
			0 & 1 
		\end{pmatrix} \right),\\
	&X_{11}(t) = \tilde{\sigma}_{11}\tilde{Y}_1+m_{11}(t)-x, \ X_{22}(t)=\tilde{\sigma}_{22}(\tilde{Y}_3-\tilde{Y}_1)+m_{22}(t)-x, \ X_{12}(t) = \tilde{\sigma}_{12}\tilde{Y}_2.
	\end{split}
\end{equation*}
Similarly to the computation of $J_{1,t}$, we can write the conditional distribution of $(X_{11}(t)X_{22}(t)-X_{12}^2(t) | X(t)=x, \tilde{Y}_3 = y)$ as
	\begin{equation*}
		\begin{split}
			&\quad \left((\tilde{\sigma}_{11}\tilde{Y}_1+m_{11}(t)-x)[\tilde{\sigma}_{22}(\tilde{Y}_3-\tilde{Y}_1)+m_{22}(t)-x] - \sigma_{12}^2\tilde{Y}_2^2\,|\, X(t)=x, \tilde{Y}_3 = y\right)\\
			&\overset{d}{=} -\tilde{\sigma}_{11}\tilde{\sigma}_{22}\Bigg[\frac{1-\tilde{\rho}}{2} \left(\mathcal{N}_1 + \frac{1}{\sqrt{2(1-\tilde{\rho})}}\left(\frac{m_{11}(t)-x}{\tilde{\sigma}_{11}}-\frac{m_{22}(t)-x}{\tilde{\sigma}_{22}}\right)\right)^2 +\tilde{\rho} \mathcal{N}_2^2- \frac{\left(y+\tilde{a}_t(x)\right)^2}{4} \Bigg],
		\end{split}
	\end{equation*}
where $\tilde{a}_t(x)$ is given in \eqref{eq:a_t}, and $\mathcal{N}_1$ and $\mathcal{N}_2$ are independent standard normal random variables. Let 
\begin{equation}\label{eq:W_tilde}
\tilde{W} = \frac{1-\tilde{\rho}}{2} \left(\mathcal{N}_1 + \frac{1}{\sqrt{2(1-\tilde{\rho})}}\left(\frac{m_{11}(t)-x}{\tilde{\sigma}_{11}}-\frac{m_{22}(t)-x}{\tilde{\sigma}_{22}}\right)\right)^2 +\tilde{\rho} \mathcal{N}_2^2.
\end{equation}
Then $\tilde{W}$ has a generalized chi-squared distribution. Applying Lemma \ref{lemma:chi-square} with $\alpha_1=\frac{1-\tilde{\rho}}{2}$, $\alpha_2=\tilde{\rho}$ and $\beta = \frac{1}{\sqrt{2(1-\tilde{\rho})}}\left(\frac{m_{11}(t)-x}{\tilde{\sigma}_{11}}-\frac{m_{22}(t)-x}{\tilde{\sigma}_{22}}\right)$, we have the pdf of $\tilde{W}$ as follows:
\begin{equation*}
	\begin{split}
		f_{\tilde{W}}(w) = \sum_{k=0}^\infty \frac{(-1)^k\tilde{c}_{k,t}(x)w^k}{2^{k+1}k!},
	\end{split}
\end{equation*}
where $\tilde{c}_{k,t}(x)$ are given in \eqref{eq:a_t}. Similarly to computing $J_{1,t}$, we obtain
\begin{equation}\label{eq:J2E}
	\begin{split}
		J_{2,t}(x) &= \frac{-\tilde{\sigma}_{11}\tilde{\sigma}_{22}}{2\sqrt{\pi(1+\tilde{\rho})}} \int_{-\infty}^{-\tilde{a}_t(x)} \E \left[ \left(\tilde{W}- \frac{1}{4}\left(y+\tilde{a}_t(x)\right)^2\right)\mathbbm{1}_{\{\tilde{W}< \frac{1}{4}\left(y+\tilde{a}_t(x)\right)^2\}}\right] e^{-\frac{y^2}{4(1+\tilde{\rho})}}dy\\
		&=\tilde{\sigma}_{11}\tilde{\sigma}_{22} \sum_{k=0}^\infty \frac{(-1)^{k+1}\tilde{a}_t^{2k+4}(x)\tilde{c}_{k,t}(x)}{2^{3k+5}(k+2)!}\sum_{j=0}^{2k+4} {2k+4 \choose j}\left(\frac{\sqrt{2(1+\tilde{\rho})}}{\tilde{a}_t(x)}\right)^jQ_j\left(\frac{-\tilde{a}_t(x)}{\sqrt{2(1+\tilde{\rho})}}\right).
	\end{split}
\end{equation}
Note that, for each $t$, $\nabla X(t)$ is independent of $X(t)$ and $\nabla^2 X(t)$ by stationarity. Therefore, by the Kac-Rice formula,
\begin{equation*}
	\begin{split}
		\E[M_u(X,T)] &= \int_T \E\left[|{\rm det} (\nabla^2 X(t))| \mathbbm{1}_{\{\nabla^2 X(t)\prec 0\}}\mathbbm{1}_{\{X(t)\ge u\}} | \nabla X(t)=(0,0)\right] p_{\nabla X(t)}(0,0) dt\\
		&= \frac{1}{(2\pi)^{3/2}}\int_Tdt\int_u^\infty e^{-\frac{(x-m(t))^2+m_1^2(t)+m_2^2(t)}{2}}J_{2,t}(x) \, dx dt.
	\end{split}
\end{equation*} 
Lastly, taking the derivative in \eqref{eq:Palm distr Euclidean}, we obtain the peak height density 
\begin{equation}\label{eq:HD-2D}
	\begin{split}
		h_t(x) = \frac{\phi(x-m(t))\E\left[|{\rm det} (\nabla^2 X(t))| \mathbbm{1}_{\{\nabla^2 X(t)\prec 0\}} | X(t)=x\right]}{\E\left[|{\rm det} (\nabla^2 X(t))|\mathbbm{1}_{\{\nabla^2 X(t)\prec 0\}}\right]}= \frac{\phi(x-m(t))J_{2,t}(x)}{J_{1,t}}.
	\end{split}
\end{equation}
	\end{proof}



\subsection{The 2D centered case}
Here, we consider centered planar Gaussian fields, where $m(t)\equiv 0$. Applying Lemma \ref{lemma:chi-square} with $\alpha_1=\frac{1-\rho}{2}$, $\alpha_2=\rho$ and $\beta=0$, by Remark \ref{remark:chi}, we obtain that the pdf of $W$ is given by $f_W(w) = \sum_{k=0}^\infty \frac{c_{k,t}(-w)^k}{2^{k+1}k!}$, where 
\begin{equation*}
c_{k,t} = \frac{\sqrt{2}}{\pi\sqrt{\rho(1-\rho)}}\sum_{i=0}^k \frac{2^i\Gamma\left(i+\frac{1}{2}\right)\Gamma\left(k-i+\frac{1}{2}\right)}{i!(k-i)!(1-\rho)^i\rho^{k-i}}.
\end{equation*}
Since $a_t=0$, it follows from \eqref{eq:J1t} that 
\begin{equation*}
	\begin{split}
J_{1,t}&=\sigma_{11}\sigma_{22}\sum_{k=0}^\infty \frac{(-1)^{k+1}c_{k,t}}{2^{3k+5}(k+2)!}  \left(\sqrt{2(1+\rho)}\right)^{2k+4}Q_{2k+4}(0)\\
&= \sigma_{11}\sigma_{22}\sum_{k=0}^\infty \frac{(-1)^{k+1}c_{k,t} (1+\rho)^{k+2} (2k+3)!!}{2^{2k+4}(k+2)!}.
\end{split}
\end{equation*}
Note that $c_{k,t}$ and $J_{1,t}$ here do not depend on $t$. In particular, by \eqref{eq:LocalMax-2D-1}, we obtain a specific expression for the expected number of local maxima $\E[M(X, T)] = \frac{{\rm Area(T)}}{2\pi}J_{1,t}$.

On the other hand, applying Lemma \ref{lemma:chi-square} with $\alpha_1=\frac{1-\tilde{\rho}}{2}$, $\alpha_2=\tilde{\rho}$ and $\beta = \frac{-x}{\sqrt{2(1-\tilde{\rho})}}\left(\frac{1}{\tilde{\sigma}_{11}}-\frac{1}{\tilde{\sigma}_{22}}\right)$, we obtain that the pdf of $\tilde{W}$ is given by $f_{\tilde{W}}(w) = \sum_{k=0}^\infty \frac{\tilde{c}_{k,t}(x)(-w)^k}{2^{k+1}k!}$,
where $\tilde{c}_{k,t}(x)$ are given in \eqref{eq:a_t}. Note that, $m(t)\equiv 0$ implies
\[
\tilde{a}_t(x)=-x\left( \frac{1}{\tilde{\sigma}_{11}}+\frac{1}{\tilde{\sigma}_{22}} \right).
\]
Here, $\tilde{c}_{k,t}(x)$ and $\tilde{a}_t(x)$ do not depend on $t$. Finally, $J_{2,t}(x)$ is given in \eqref{eq:J2t}, leading to $\E[M_u(X,J)]$ and $h_t(x)$ as provided in \eqref{eq:LocalMax-2D-2} and \eqref{eq:phd1}, respectively.

\subsection{The 2D isotropic case}
We now explore the scenario where $X$ is a planar isotropic Gaussian random field. Specifically, in Corollary \ref{cor:Isotropic}, we present new results for isotropic Gaussian fields with general mean functions. Additionally, we confirm that our method, utilizing the chi-squared density, yields the same formula as that obtained through the random matrix technique \cite{CS18}.
\subsubsection{Noncentered isotropic Gaussian fields}
Due to isotropy, the covariance function of $X$ can be written as ${\rm Cov}(X(t), X(s))=\varphi(\|t-s\|^2)$ for an appropriate function $\varphi(\cdot): [0,\infty) \rightarrow \R$. We denote
\begin{equation}\label{Eq:kappa}
	\varphi'=\varphi'(0), \quad \varphi''=\varphi''(0),  \quad \kappa=-\varphi'/\sqrt{\varphi''}.
\end{equation}
It can be verified that (cf. \cite{CS18}):
\begin{equation*}
	\begin{split}
{\rm Cov}(\nabla X(t)) = -{\rm diag}(2\varphi', 2\varphi'), \quad 		{\rm Cov}(X_{11}(t), X_{22}(t), X_{12}(t)) = \begin{pmatrix} 
			12\varphi'' & 4\varphi'' & 0 \\
			4\varphi'' & 12\varphi'' & 0 \\
			0 & 0 & 4\varphi'' \\
		\end{pmatrix}.
	\end{split}
\end{equation*}
Under our assumption that ${\rm Cov}(\nabla X(t)) = I_2$, we find $\varphi'=-1/2$. Moreover, we have $\sigma_{11}^2=\sigma_{22}^2=12\varphi''$ and $\sigma_{12}^2=4\varphi''$.  It is well-known that, for an isotropic Gaussian field, the Hessian matrix $\nabla^2 X(t)$ is orthogonally invariant; that is, the distribution of $\nabla^2 X(t)$ is the same as that of $U(\nabla^2 X(t))U^T$ for any $2\times 2$ orthogonal matrix $U$.

The following result demonstrates that, for an isotropic Gaussian field, we can compute the expected number and height distribution of local maxima when incorporating a general mean function, without requiring $\nabla^2 m(t)$ to be a diagonal matrix.

\begin{corollary}\label{cor:Isotropic}
	Let $\{X(t), t\in \R^2\}$ be a smooth, unit-variance, isotropic Gaussian field with covariance ${\rm Cov}(X(t), X(s))=\varphi(\|t-s\|^2)$, mean function $m(t)$, and satisfying ${\rm Cov}(\nabla X(t)) = I_2$ (i.e., $\varphi'=-1/2$). Then, the results in Theorem \ref{theorem:2D} hold with $\sigma_{11}^2=\sigma_{22}^2=12\varphi''$, $\sigma_{12}^2=4\varphi''$, and by replacing $(m_{11}(t), m_{22}(t))$ with $(\theta_{1,t}, \theta_{2,t})$, where $\theta_{1,t}$ and $\theta_{2,t}$ are the eigenvalues of $\nabla^2 m(t)$. 
\end{corollary}

\begin{proof}
	Let $U_t$ be the orthogonal matrix such that $U_t(\nabla^2 m(t)) U_t^T = {\rm diag} (\theta_{1,t}, \theta_{2,t})$. By the orthogonally invariant property of the Hessian, we have
	\begin{equation*}
		\begin{split}
			J_{1,t} &= \E\left[|{\rm det} (\nabla^2 Z(t) + \nabla^2 m(t))| \mathbbm{1}_{\{\nabla^2 Z(t) + \nabla^2 m(t)\prec 0\}}\right]\\ 
			&=\E\left[|{\rm det} (U_t(\nabla^2 Z(t)) U_t^T+ U_t(\nabla^2 m(t))U_t^T)| \mathbbm{1}_{\{U_t(\nabla^2 Z(t)) U_t^T+ U_t(\nabla^2 m(t))U_t^T\prec 0\}}\right]\\ 
			&=\E\left[|{\rm det} (U_t(\nabla^2 Z(t)) U_t^T+ {\rm diag} (\theta_{1,t}, \theta_{2,t}))| \mathbbm{1}_{\{U_t(\nabla^2 Z(t)) U_t^T+ {\rm diag} (\theta_{1,t}, \theta_{2,t})\prec 0\}}\right]\\ 
			&=\E\left[|{\rm det} (\nabla^2 Z(t)+ {\rm diag} (\theta_{1,t}, \theta_{2,t}))| \mathbbm{1}_{\{\nabla^2 Z(t)+ {\rm diag} (\theta_{1,t}, \theta_{2,t})\prec 0\}}\right].
		\end{split}
	\end{equation*}
	Applying the same arguments to $J_{2,t}$, we see that the results in Theorem \ref{theorem:2D} hold if $m_{11}(t)$ and $m_{22}(t)$ are replaced with $\theta_{1,t}$ and $\theta_{2,t}$ respectively. 
\end{proof}

\subsubsection{Centered isotropic Gaussian fields}

Here, we assume $m(t)\equiv 0$ and use the derived method to find explicit formulas for the expected number and height distribution of local maxima. It will be shown that the derived formulas are the same as those in \cite{CS18}. Note that $\varphi'=-1/2$,  $\sigma_{11}^2=\sigma_{22}^2=12\varphi''$ and $\sigma_{12}^2=4\varphi''$. Thus
\[
\kappa=\frac{1}{2\sqrt{\varphi''}}, \quad \tilde{\sigma}_{11}^2=\tilde{\sigma}_{22}^2=12\varphi''-1,\quad \rho = \frac{1}{3}, \quad \tilde{\rho} = \frac{4\varphi''-1}{12\varphi''-1}, \quad \frac{1-\tilde{\rho}}{2} = \frac{4\varphi''}{12\varphi''-1}.
\]
This implies that $W$ defined in \eqref{eq:W} has the distribution of $\frac{1}{3}\chi_2^2$. Plugging this into \eqref{eq:J1E}, together with $a_t=0$, we obtain
\begin{equation}\label{eq:J1-Iso}
	\begin{split}
		J_{1,t} = -\frac{12\varphi''}{2\sqrt{4\pi/3}}\int_{-\infty}^0 \frac{1}{3}\E \left[ \left(\chi_2^2 - \frac{3y^2}{4}\right) \mathbbm{1}_{\{\chi_2^2 < \frac{3y^2}{4}\}} \right]  e^{-\frac{3y^2}{16}}dy= \frac{1}{\sqrt{3}\kappa^2},
	\end{split}
\end{equation}
where we have used the fact that, for a constant $c$,
\begin{equation}\label{eq:chi-expectation}
\E \left[ \left(\chi_2^2- c\right)\mathbbm{1}_{\{\chi_2^2< c\}} \right] = \int_0^c (w-c)\frac{1}{2}e^{-\frac{w}{2}}dw =2-c -2e^{-\frac{c}{2}}.
\end{equation}

Similarly, we see that $\tilde{W}$ defined in \eqref{eq:W_tilde} has the distribution of $\frac{4\varphi''}{12\varphi''-1} \chi_2^2$. Plugging this into \eqref{eq:J2E}, together with $\tilde{a}_t(x)=-\frac{2x}{\sqrt{12\varphi''-1}}$ and \eqref{eq:chi-expectation}, we obtain
\begin{equation}\label{eq:J2-expand}
	\begin{split}
		J_{2,t}(x)
		&=  \frac{-2\varphi''\sqrt{12\varphi''-1}}{\sqrt{2\pi(8\varphi''-1)}} \int_{-\infty}^{\frac{2x}{\sqrt{12\varphi''-1}}}\Bigg[2-\frac{12\varphi''-1}{16\varphi''}\left(y-\frac{2x}{\sqrt{12\varphi''-1}}\right)^2\\
		&\qquad -2e^{-\frac{12\varphi''-1}{32\varphi''}\left(y-\frac{2x}{\sqrt{12\varphi''-1}}\right)^2} \Bigg]e^{-\frac{y^2(12\varphi''-1)}{8(8\varphi''-1)}}dy.
	\end{split}
\end{equation}
Note that 
\begin{equation*}
	\begin{split}
\int_{-\infty}^{\frac{2x}{\sqrt{12\varphi''-1}}}e^{-\frac{y^2(12\varphi''-1)}{8(8\varphi''-1)}}dy &=  2\sqrt{2\pi}\sqrt{\frac{8\varphi''-1}{12\varphi''-1}}\Phi\left(\frac{x}{\sqrt{8\varphi''-1}}\right),\\
\int_{-\infty}^{\frac{2x}{\sqrt{12\varphi''-1}}}ye^{-\frac{y^2(12\varphi''-1)}{8(8\varphi''-1)}}dy &=  -\frac{4(8\varphi''-1)}{12\varphi''-1}e^{-\frac{x^2}{2(8\varphi''-1)}},\\
\int_{-\infty}^{\frac{2x}{\sqrt{12\varphi''-1}}}y^2e^{-\frac{y^2(12\varphi''-1)}{8(8\varphi''-1)}}dy &=  -  \frac{8x(8\varphi''-1)}{(12\varphi''-1)^{3/2}}e^{-\frac{x^2}{2(8\varphi''-1)}}  + 8\sqrt{2\pi}\left(\frac{8\varphi''-1}{12\varphi''-1}\right)^{3/2}\Phi\left(\frac{x}{\sqrt{8\varphi''-1}}\right),
\end{split}
\end{equation*}
and 
\begin{equation*}
	\begin{split}
		&\quad \int_{-\infty}^{\frac{2x}{\sqrt{12\varphi''-1}}}e^{-\frac{12\varphi''-1}{32\varphi''}\left(y-\frac{2x}{\sqrt{12\varphi''-1}}\right)^2} e^{-\frac{y^2(12\varphi''-1)}{8(8\varphi''-1)}}dy \\
		&= 4\sqrt{2\pi}  \frac{\sqrt{\varphi''(8\varphi''-1)}}{12\varphi''-1} e^{-\frac{x^2}{2(12\varphi''-1)}}\Phi\left( \sqrt{\frac{4\varphi''}{(12\varphi''-1)(8\varphi''-1)}}x \right).
	\end{split}
\end{equation*}
Plugging these into \eqref{eq:J2-expand}, we obtain
\begin{equation*}
	\begin{split}
		J_{2,t}(x)
		&=   (x^2-1)\Phi\left(\frac{x}{\sqrt{8\varphi''-1}}\right) + \frac{\sqrt{8\varphi''-1}}{\sqrt{2\pi}} xe^{-\frac{x^2}{2(8\varphi''-1)}} \\
		&\quad + 16  \frac{\varphi''\sqrt{\varphi''}}{\sqrt{12\varphi''-1}} e^{-\frac{x^2}{2(12\varphi''-1)}}\Phi\left( \sqrt{\frac{4\varphi''}{(12\varphi''-1)(8\varphi''-1)}}x \right).
	\end{split}
\end{equation*}
Since $\varphi''=\frac{1}{4\kappa^2}$, we obtain
\begin{equation}\label{eq:J2-Iso}
	\begin{split}
		J_{2,t}(x)&=  (x^2-1)\Phi\left(\frac{\kappa x}{\sqrt{2-\kappa^2}}\right) + \frac{\sqrt{2-\kappa^2}}{\sqrt{2\pi}\kappa } xe^{-\frac{\kappa^2x^2}{2(2-\kappa^2)}} \\
		&\quad +  \frac{2}{\kappa^2\sqrt{3-\kappa^2}} e^{-\frac{\kappa^2x^2}{2(3-\kappa^2)}}\Phi\left( \frac{\kappa x}{\sqrt{(3-\kappa^2)(2-\kappa^2)}} \right).
	\end{split}
\end{equation}
By \eqref{eq:phd1}, taking the ratio between \eqref{eq:J2-Iso} and \eqref{eq:J1-Iso}, we obtain the peak height density
\begin{equation*}
	\begin{split}
		h_t(x)
		&= \sqrt{3}\kappa^2(x^2-1)\phi(x)\Phi\left(\frac{\kappa x}{\sqrt{2-\kappa^2}} \right) + \frac{\kappa x\sqrt{3(2-\kappa^2)}}{2\pi}e^{-\frac{x^2}{2-\kappa^2}} \\
		&\quad +\frac{\sqrt{6}}{\sqrt{\pi(3-\kappa^2)}}e^{-\frac{3x^2}{2(3-\kappa^2)}}\Phi\left(\frac{\kappa x}{\sqrt{(3-\kappa^2)(2-\kappa^2)}} \right),
	\end{split}
\end{equation*}
which is consistent with that derived in Example 3.8 in \cite{CS18} using random matrix techniques.

	\section*{Acknowledgments}
	The author was supported by NSF Grants DMS-1902432 and DMS-2220523, and Simons Foundation Collaboration Grant \#854127. He is grateful to Armin Schwartzman from UC San Diego and Xuan Yang for their helpful discussions and suggestions. He also thanks the two anonymous referees for their valuable comments that led to improvements in the paper.
	
	\begin{small}
		
	\end{small}
	
	\bigskip
	
	\begin{quote}
		\begin{small}
			
			\textsc{Dan Cheng}\\
			School of Mathematical and Statistical Sciences \\
			Arizona State University\\
			900 S. Palm Walk\\
			Tempe, AZ 85281, U.S.A.\\
			E-mail: \texttt{cheng.stats@gmail.com}

		\end{small}
	\end{quote}
	
\end{document}